\documentclass[11pt]{article}
\usepackage{amsmath}
\usepackage{latexsym}
\usepackage{graphicx}
\usepackage{amssymb}
\newtheorem{theorem}{Theorem}

\newtheorem{question}{Question}
\newtheorem{lemma}{Lemma}
\newtheorem{observation}{Observation}
\newtheorem{corollary}{Corollary}
\author{Eli Berger}
\title{Dynamic Monopolies of Constant Size}
\begin{document}
\maketitle

\begin{abstract}

The paper deals with a polling game on a graph. Initially,
each vertex is colored white or black. At each round, each
vertex is colored 
by the color shared by the majority of vertices in its neighborhood,
at the previous round.
(All recolorings are done simultaneously).
We say that a set $W_0$ of
vertices is a {\em dynamic monopoly} or {\em dynamo} if 
starting the game with the vertices of $W_0$ colored 
white, the entire system is white after a finite number of
rounds. Peleg \cite{Peleg} asked how small a dynamic monopoly  
may be as a function of the number of vertices. We show that
the answer is O(1). 

\end {abstract}

\section{Introduction}

Let $G=(V,E)$ be a simple undirected graph
and $W_0$ a subset of $V$. Consider the following repetitive 
polling game. At round 0 the vertices of $W_0$ are colored white
and the other vertices are colored black.
At each round, each vertex $v$ is colored 
according to the following rule.
If at round $r$ the vertex $v$ has
more than half of its neighbors colored $c$, then 
at round $r+1$ the vertex $v$ will be colored $c$.
If at round $r$ the vertex $v$ has
exactly half of its neighbors colored white 
and half of its neighbors colored black, then we say there is
a tie. In this case $v$ is colored at round $r+1$ by the same
color it had at round $r$. (Peleg considered other models for
dealing with ties. We will refer to these models in section \ref{remarks}.
Additional models and further study of this game may be found at 
\cite{GO80}, \cite{PS83}, \cite{LPRS93}, \cite{BP95} and \cite{BBPP96}.)
If there exists a finite $r$ so that at round $r$ all vertices in $V$
are white, then we say that $W_0$ is a {\em dynamic monopoly},
abbreviated {\em dynamo}.

In this paper we prove

\begin{theorem}\label{main}
For every natural number $n$ there exists a graph with more than
$n$ vertices and with a dynamic monopoly of 18 vertices.
\end{theorem}

We shall use the following notation:
If $v \in V$ then $N(v)$ denotes the set of neighbors of $v$. We call
$d(v) = |N(v)|$ the {\em degree} of $v$. For every $r = 0,1 \ldots$
we define $C_r$ as a function from $V$ to $\{ {\cal B}, {\cal W} \}$, so
that $C_r(v) = {\cal W}$ if $v$ is white at round $r$ and
$C_r(v) = {\cal B}$ if $v$ is black at this round.
We also define $W_r = C_r^{-1}({\cal W})$,  $B_r = C_r^{-1}({\cal B})$,
$T_r = W_r \cap W_{r-1} (r>0)$ and $S_r = T_1 \cup \ldots \cup T_r$

\section{Proof of Theorem \ref{main}}\label{proofmain}

\begin{figure}[tbhp]
\begin{center}\leavevmode
\includegraphics [width=\linewidth] {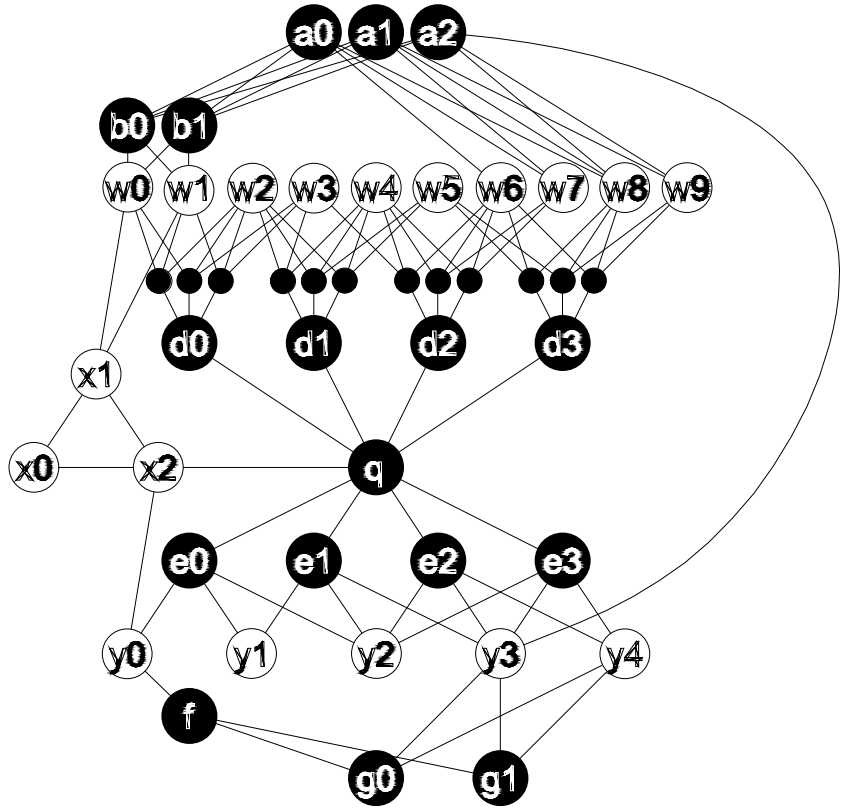}
\caption{The graph $J$. (The small black circles are the vertices $c_0 \ldots c_{11}$.)}
\end{center}\label{JJ}
\end{figure}

Let $J=(V_J,E_J)$ be the graph in figure 1. Let 

$$ W_0 = \{ w_0, \ldots, w_9, x_0, \ldots, x_2, y_0, \ldots, y_4  \} $$

and let $U = W_0 \cup \{ q \}$ and $D = V_J - U$. We construct 
a graph $J_n$ by duplicating $n$ times the vertices in $D$.
That is,

$$ J_n = (V_n,E_n) $$

where

$$ V_n = U \cup [n] \times D $$

and

\begin{eqnarray} 
E_n = \{ (u,v) \in J : u,v \in U \} \cup
\{ (u,(i,v)) : (u,v) \in J , u \in U , v \in D , i \in [n] \}
\nonumber
\\ \cup \{ ((i,u),(i,v)) : (u,v) \in J , \: u,v \in D , i \in [n] \}
\nonumber 
\end{eqnarray}

(Here, as usual, $[n]$ denotes the set $\{ 1 \ldots n \}$).

Note that for reasons of symmetry, at a given round,
all copies of a vertex in $J$ have
the same color. Thus we may write
``$y_0$ is white at round 3'' instead of
``$(i,y_0)$ is white at round 3 for every $i \in [n]$'' etc.

The following table describes the evolution of $J_n$. 
The symbol 1 stands for white and 0 stands for black.
Note that the table does {\em not} depend on $n$. 
(This property is peculiar to the graph $J$.
In general graphs duplication of vertices may
change the pattern of evolution of the graph).

$
\begin{array}{ccccccccccc}
r & a_{012} & b_{01} & c_0 \ldots c_{11} & d_{0123} &
e_{0123} & f & g_{01} & q & w_0 \ldots w_9 &
y_{01234}  \\
0&000&00&000000000000&0000&0000&0&00&0&1111111111&11111\\
1&111&00&111111111111&0000&1111&0&11&0&0000000000&00000\\                         
2&000&11&000000000000&1111&0000&1&00&1&1111111111&11111\\                         
3&111&00&111111111111&0000&1111&0&11&1&1100000000&10000\\                         
4&000&11&100000000000&1111&1000&1&00&1&1111111111&11111\\                         
5&111&00&111111111111&1000&1111&0&11&1&1100000000&11000\\                         
6&000&11&111000000000&1111&1100&1&00&1&1111111111&11111\\                         
7&111&00&111111111111&1000&1111&0&11&1&1111000000&11100\\                         
8&000&11&111100000000&1111&1111&1&00&1&1111111111&11111\\                         
9&111&00&111111111111&1100&1111&0&11&1&1111000000&11111\\                         
10&000&11&111111000000&1111&1111&1&11&1&1111111111&11111\\                         
11&111&00&111111111111&1100&1111&1&11&1&1111110000&11111\\                         
12&000&11&111111100000&1111&1111&1&11&1&1111111111&11111\\                         
13&111&00&111111111111&1110&1111&1&11&1&1111110000&11111\\                         
14&000&11&111111111000&1111&1111&1&11&1&1111111111&11111\\                         
15&111&00&111111111111&1110&1111&1&11&1&1111111100&11111\\                         
16&000&11&111111111100&1111&1111&1&11&1&1111111111&11111\\                         
17&111&00&111111111111&1111&1111&1&11&1&1111111100&11111\\                         
18&000&11&111111111111&1111&1111&1&11&1&1111111111&11111\\                         
19&111&00&111111111111&1111&1111&1&11&1&1111111111&11111\\                         
20&111&11&111111111111&1111&1111&1&11&1&1111111111&11111\\                         
21&111&11&111111111111&1111&1111&1&11&1&1111111111&11111
\end{array}
$

The table shows that at round 20 the entire system is white
and therefore $W_0$ is a dynamo. 
The reader may go through the table by himself, but in order
to facilitate the understanding of what happens in the table
let us add some explanations as to
the mechanism of ``conquest'' used in this graph.  

We say that round $j$ { \em dominates } round $i$ if
$W_i \subseteq W_j$.

We shall make use of the following obvious fact:

\begin{observation}\label{monotony}
If round $j$ dominates round $i$ ($i,j = 0,1 \ldots$) then
round $j+1$ dominates round $i+1$.
\end{observation} 



 




By applying this observation $k$ times, we find that
if round $j$ dominates round $i$ then
round $j+k$ dominates round $i+k$ ($i,j,k = 0,1 \ldots $).
By looking at the table one can see that in the graph $J_n$ round
2 dominates round 0 and thus we have 

\begin{corollary}\label{blinking}
Round $k+2$ dominates round $k$ in $J_n$ for every $k = 0,1 \ldots$ 
\end{corollary} 

We say that a vertex $v$ {\em blinks} at round $r$ if $C_{r+2i}(v) = {\cal W}$
 for every $i=0,1 \ldots$.
We say that a vertex $v$ is {\em conquered} at round $r$ if $C_{r+i}(v) = {\cal W}$
 for every $i=0,1 \ldots$. 
Examining rounds $0$ to $3$ in the table and using Corollary
\ref{blinking} one can see that $x_0, x_1$ and $x_2$ are
conquered at round 0, and in addition $q, w_0, w_1$ and $y_0$
are conquered at round 2. Furthermore, every vertex in
$J_n$ blinks either at round 1 or at round 2.

Finally, we have

\begin{lemma}\label{tieconquer}
If at round $r$ a vertex $v$ in $J_n$ has at least half of
its neighbors conquered then
$v$ is conquered at round $r+2$.
\end{lemma} 

{\em Proof:}
Every vertex in
$J_n$ blinks either at round 1 or at round 2,
and hence $v$ is white either at round $r+1$ or at round $r+2$.
From this round on, at least half of the neighbors of $v$ 
are white, so $v$ will stay white.

$\Box$

Now the vertices will be conquered in the following order:

$x_0,x_1,x_2,q,w_0,w_1,y_0,c_0,e_0,d_0,y_1,c_1,c_2,e_1,w_2
w_3,y_2,c_3,e_2,e_3,d_1,y_3,y_4$, 

$c_4,c_5,g_0,g_1,f,w_4,w_5,c_6,d_2,c_7,c_8,w_6,w_7,c_9,d_3,c_{10},c_{11},w_8,
w_9$,

$a_0,a_1,a_2,b_0,b_1$.

Eventually, the entire graph is colored white. $J_n$ is a graph
with $19+27n > n$ vertices and $W_0$ is a dynamo of size 18,
proving Theorem \ref{main}.

\section{Questions and Remarks}\label{remarks}

The result of Section \ref{proofmain} gives rise to
the following questions:


\begin{question}
Does there exist an infinite graph with a finite dynamo?
\end{question}

The answer is {\em no}. 
This follows from the following theorem:
\begin{theorem}
If $W_0$
is finite then 
$T_r$ is finite for all $r=1,2 \ldots$. 
Moreover, every vertex in $T_r$ has a finite degree.
\end{theorem}

{\em Proof:}
The proof is by induction on $r$. For $r=1$ the theorem is true 
because every vertex $v \in W_0$ with an infinite degree becomes 
black at round 1.
For $r>1$, if $C_{r-1}(v) = {\cal W}$ and $v$ has 
an infinite degree $\lambda$
then by the induction hypotheses $C_{r-2}(v) = {\cal B}$ and 
$|N(v) \cap B_{r-2}| < \lambda$.
Hence $|N(v) \cap W_{r-1}| \leq |N(v) \cap B_{r-2}| + |T_{r-1}| < \lambda$
and $C_r(v) = {\cal B}$. 

If $v \in T_r$ has a finite degree then $v$ has a neighbor 
in $T_{r-1}$. 
By the induction hypotheses only finitely many vertices have
such a neighbor, and thus $T_r$ is finite.

$\Box$ 

The next question deals with other models considered by Peleg:

\begin{question}\label{2G}
Do we still have a dynamo of size O(1) if we change the rules
of dealing with ties? (e.g. if a vertex becomes black whenever
there is a tie.)
\end{question}

The answer here is {\em yes}. If $G=(V,E)$ is a graph, 
introduce a new vertex $v'$ for every $v \in V$ and consider
the graph $\hat{G} = (\hat{V},\hat{E})$ where

$$\hat{V} = \{v,v': v \in V\}$$ 

and

$$ \hat{E} = E \cup \{(u',v'):(u,v) \in E \} \cup \{(v,v'): 2|d(v) \} $$

If $W_0$ is a dynamo of $G$ according to the
model in Theorem \ref{main}, then it is easy to prove that 
$\hat{W_0} = \{v,v':v \in W_0\}$
is a dynamo of $\hat{G}$. But all vertices of $\hat{G}$ have odd degrees, and
thus ties are not possible and $\hat{W_0}$
is a dynamo of $\hat{G}$ according to {\em any} rule 
of dealing with ties.

Therefore, for every $n = 1,2 \ldots $ the graph $\hat{J_n}$ has a dynamo
of size 36.



\section{Another Model}

Let $\rho>1$ be a real number. Consider the following 
model, which will henceforth be called {\em the $\rho$-model}. At 
every round, for every vertex $v$ with $b$ neighbors colored black
and $w$ neighbors colored white, if $w>\rho b$ then $v$ is colored white
at the next round, otherwise it is black. For the sake of simplicity 
we will assume that $\rho$ is irrational and that there are no
isolated vertices, so that $w=\rho b$ is impossble.


The most interesting question regarding this model is whether there exist
graphs with O(1) dynamo like in Theorem \ref{main}. This question is as
yet open. We only have some partial results, which can be summarized as
follows:

 \begin{enumerate}
 \renewcommand{\theenumi}{\roman{enumi}}
 
 \item If $\rho$ is big enough then the size of a dynamo
 is $\Omega(\sqrt{n})$.
 
 \item If $\rho$ is small enough then 
 there exist graphs in which
 the size of a dynamo
 is $O(\log n)$.
 
 \item If there exist graphs with O(1) dynamo then the number of 
 rounds needed until the entire system becomes white is
 $\Omega(\log n)$.
  
 \end{enumerate}

More explicitly:

\begin{theorem}\label{rho3}
Let $ \rho > 3$. If a graph with $n$ vertices has a dynamo of size $k$
in the $\rho$-model then 
$$n < k^2 $$
\end{theorem}

{\em proof:}

For every $r = 1,2, \ldots$, let $(S_r,\overline S_r)$ be 
the set of edges with one vertex in $S_r$ and the other
not in $S_r$. Call $s_r = |S_r|+|(S_r,\overline S_r)|$. Note that
$S_1$ is the set of vertices
which are white at both round 0 and round 1. Every $v \in S_1$
is connected to at most $k-|S_1|$ vertices in $W_0 \setminus S_1$
and at most $\frac{k-1}{\rho}<k-1$ vertices outside of $W_0$.
Therefore we have  
$$s_1 < |S_1| + |S_1|(k-|S_1|+k-1) = k^2 - (k-|S_1|)^2 \leq k^2$$ 

Thus all we need is to show $s_{r+1} \leq s_r$ and we are done.

Let $r$ be fixed. By definition $S_r \subseteq S_{r+1}$. 
Let $\Delta = S_{r+1} \setminus S_r$, and let $v \in \Delta$.
More than $\frac{3}{4}$ of the neighbors of $v$ are white 
at round $r$ and more than $\frac{3}{4}$ of the neighbors of $v$ are white 
at round $r-1$. Thus more than $\frac{1}{2}$ of the neighbors of $v$
belong to $S_r$. We therefore have 
$$|(S_r,\overline S_r) \setminus (S_{r+1},\overline S_{r+1})| -
 |(S_{r+1},\overline S_{r+1}) \setminus (S_r,\overline S_r)| 
\geq |\Delta| $$
which implies $s_{r+1} \leq s_r$. By induction $s_r < k^2$
for all $r$. If we begin with a dynamo then for some finite $m$
we have $S_m = V$ and $n = s_m < k^2$
$\Box$

\begin{theorem}
Let $\rho>1$. If $|W_0|=k$ and $W_m=V$ (the set of all vertices), 
then the number $e$ of edges in the graph satisfies
$$e< k^2 (\frac{2\rho}{\rho - 1})^m $$
\end{theorem}

{\em proof:}

Let $d_r$ denote the sum of the degrees of the vertices in $S_r$.
Recall that every $v \in S_1$ is white at both round 0 and round 1, and thus
$|N(v) \cap B_0| < k$ and $d(v)<k$. Therefore, 
$d_1 < 2k^2$. Again, let $r$ be fixed, let $\Delta$ be as 
in the proof of Theorem \ref{rho3} and let $v \in \Delta$.
More than $\frac{\rho}{\rho+1}$ of the neighbors of $v$ are white 
at round $r$ and more than $\frac{\rho}{\rho+1}$ 
of the neighbors of $v$ are white 
at round $r-1$. Thus more than $\frac{\rho-1}{rho+1}$ 
of the neighbors of $v$ belong to $S_r$. Therefore,
we have

$$ d_{r+1} < d_r+\frac{\rho+1}{\rho-1}d_r =  \frac{2\rho}{\rho-1}d_r $$

By induction $d_r < 2k^2 (\frac{2\rho}{\rho - 1})^{r-1}$. 
If the entire system is white at round $m$ then $d_{m+1} = 2e$
and thus we have 

$$e< k^2 (\frac{2\rho}{\rho - 1})^m $$

$\Box$

\begin{theorem}\label{small_rho}
Let $1<\rho<\frac{257}{256}$. 
For every integer $n>5$ there exists in the $\rho$ model
a graph with more than
$2^n$ vertices and with a dynamo of size $30(n-5)+36$.
\end{theorem}

{\em Outline of proof:}

Let $\hat{J}$ be as defined in the answer to Question \ref{2G}.
Construct $\tilde{J}$ by eliminating $f$ from $\hat{J}$ and connecting
$f'$ to $y_0$ and $g_1$ (but {\em not} to $g_0$). Note that
in $\tilde{J}$ the vertex $g_0$ is connected only to $y_3$ and to $y_4$.

In figure 2, the upper graph is a part of $\hat{J}$. 
The lower graph is the corresponding part in $\tilde{J}$. 
The rest of $\tilde{J}$ is identical to the rest of $\hat{J}$.

\begin{figure}[tbhp]\label{diet}
\begin{center}\leavevmode
\includegraphics[width=\linewidth]{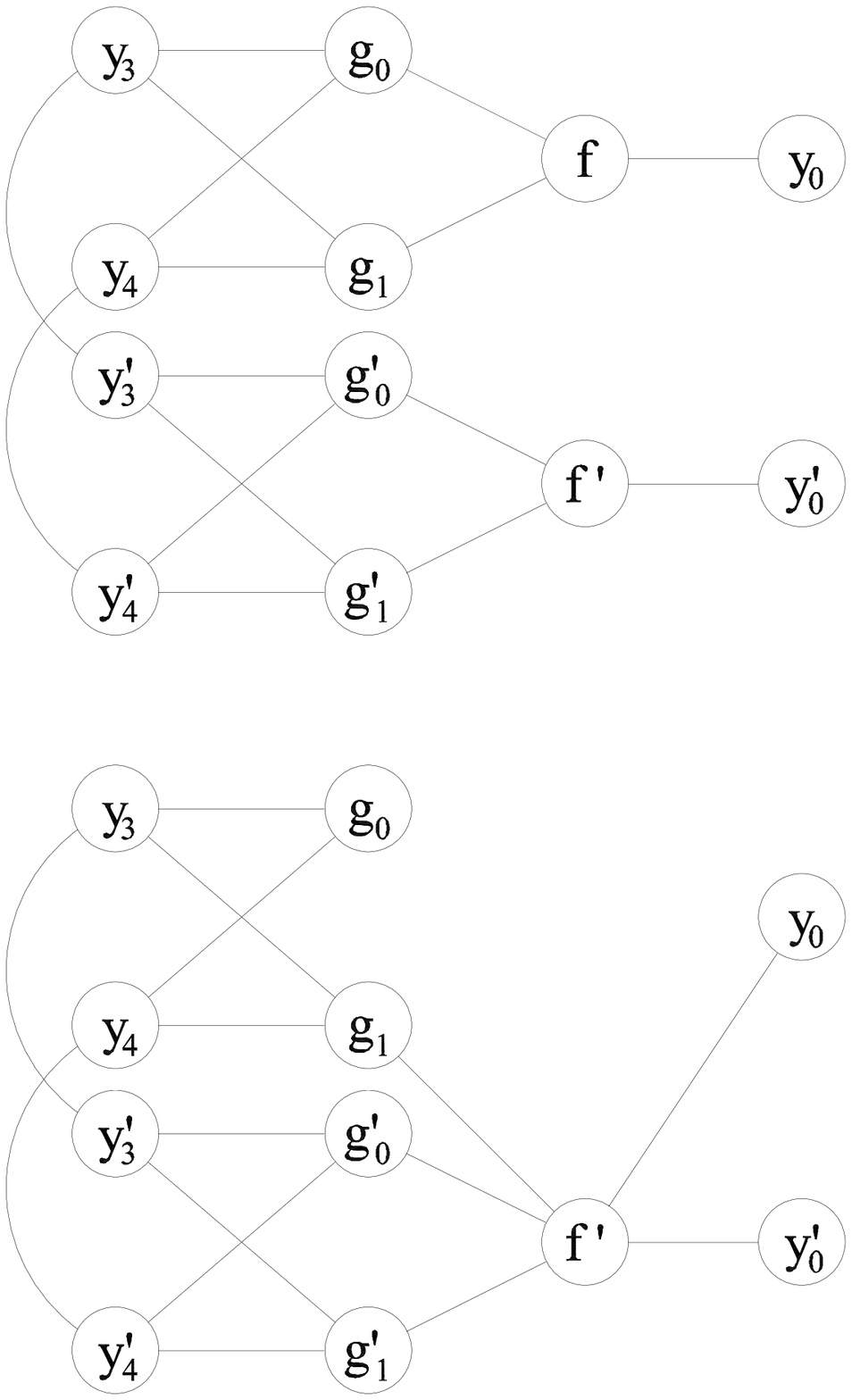}
\end{center}
\end{figure}

Construct $\tilde{J}_{32}, \tilde{J}_{64}, \ldots \tilde{J}_{2^n}$ 
as in the construction of $J_n$, where
the duplicated vertices are all black vertices except for $q$ and $q'$.
(Note that the graphs are constructed separately,
namely, the sets of vertices of $\tilde{J}_{2^i}$ and 
$\tilde{J}_{2^j}$ are disjoint  
for $i \neq j$.) Now connect the graphs in the following way.
First, eliminate the copies of $x_0,x_1,x_2$ from all graphs
except for $\tilde{J}_{32}$. 
Note that in $\tilde{J}_{2^i}$ there are $2^i$ copies
of $g_0$ (when $i = 5, \ldots n-1$).
Divide them into 32 disjoint sets $P_0, \ldots P_{31}$, of 
size $2^{i-5}$ each.
Now connect the vertices in $P_0$ to the copy of $q$ in
$\tilde{J}_{2^{i+1}}$, connect $P_1$ to the copy of $q'$,
and connect each one of $P_2 \ldots P_{31}$ to a respective
white vertex in $\tilde{J}_{2^{i+1}}$ (see in figure 3).

\begin{figure}[tbhp]\label{bigGraph}
\begin{center}\leavevmode
\includegraphics[width=\linewidth]{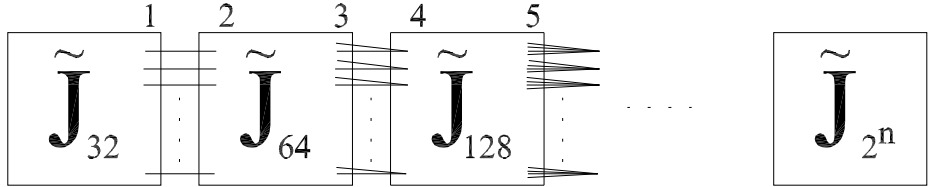}
\caption{This figure illustrates the graph used in the proof of 
Theorm \ref{small_rho}. The vertices under the numeral 1 are the 32
copies of $g_0$ in $\tilde{J}_{32}$. Under the numeral 2 are the 32
unduplicated vertices in $\tilde{J}_{64}$ ($q$, $q'$ and the initiallly
white vertices). Under the numeral 3 are the 64 copies of $g_0$ in $\tilde{J}_{64}$,
under the numeral 4 are the 32 unduplicated vertices in $\tilde{J}_{128}$,
under the numeral 5 are the 128 copies of $g_0$ in $\tilde{J}_{128}$,
and so on.}
\end{center}
\end{figure}

It is possible to verify the following:

 \begin{enumerate}
 \renewcommand{\theenumi}{\roman{enumi}}
 
 \item All vertices of the obtained graph blink either at
 round 1 or at round 2.
 
 \item All vertices of $K_{32}$ are eventually conquered. (The 
 evolution of this conquest is similar to the one in
 Theorem \ref{main}.)

 \item If all copies of $g_0$ in $\tilde{J}_{2^i}$ are conquered at a 
 certain round, then all vertices of $\tilde{J}_{2^{i+1}}$ are
 eventually conquered. (Again, the evolution is similar to the 
 one in Theorem \ref{main}. Note that we need the bound 
 $ \rho<\frac{257}{256} $ in order to have $q$ and $q'$ conquered.)

 \end{enumerate}

Thus all vertices are eventually conquered.
The theorem follows upon noticing that our graph has
more than $2^n$ vertices, and the size of the dynamo is $30(n-5)+36$.
$\Box$

{\em Acknowledgement:}

I would like to thank Ron Aharoni and Ron Holzman for helping me with the representation.


\begin{thebibliography}{9}

\bibitem{Peleg}D. Peleg, Size bounds for dynamic monopolies, {\em Discrete Applied Mathematics}, Vol: 86, Issue: 2-3, 
September 1998 (262-273).

\bibitem{GO80} E. Goles and J. Olivos, Periodic behavior of generalized threshold functions, {\em Discrete Applied Mathematics},
30:187-189, 1980.

\bibitem{PS83}S. Poljak and M. Sura, On periodic behavior in societies with symmetric influences, {\em Combinatorica},
3:119-121, 1983.

\bibitem{LPRS93} N. Linial, D. Peleg, Y. Rabinovich, and M. Saks, Sphere packing and local
majorities in graphs. In 2nd ISTCS, pages 141-149, IEEE Computer Soc. Press, June 1993.

\bibitem{BP95}J-C. Bermond and D. Peleg, The power of small coalitions in graphs, {\em Proc. 2nd
Colloc. on Structural Information and Communication Complexity}, Olympia, Greece, June 1995, 
Carleton Univ. Press, 173-184.

\bibitem{BBPP96}J-C. Bermond, J. Bond, D. Peleg, and S. Perennes, Tight bounds on the size of 
2-monopolies, {\em Proc. 3rd Colloc. on Structural Information and Communication Complexity},
June 1996, Siena, Italy. 

\end{thebibliography}
\end{document}